\documentstyle[amssymb]{article}
\date{}
\begin{document}
{\noindent\Large \texttt{Boolean $2$-designs and the\\ embedding
of a $2$-design in a group} \footnote{AMS MCS 05B05, 05B25.}}

\bigskip\noindent
{{Andrea Caggegi}\footnote{Supported by M.I.U.R.}\\
\footnotesize{Dipartimento di Metodi e Modelli Matematici}\vspace*{-1,5mm}\\
\footnotesize{Viale delle Scienze Ed. 8, I-90128 Palermo (Italy)}\vspace*{-1,5mm}\\
\footnotesize{caggegi@unipa.it}}

\medskip\noindent
{{Alfonso Di Bartolo}\footnote{Supported by Universit\`a di Palermo (Co.R.I.).}\\
\footnotesize{Dipartimento di Matematica e Applicazioni}\vspace*{-1,5mm}\\
\footnotesize{Via Archirafi 34, I-90123 Palermo (Italy)}\vspace*{-1,5mm}\\
\footnotesize{alfonso@math.unipa.it}}

\medskip\noindent
{{Giovanni Falcone}\footnote{Supported by M.I.U.R., Universit\`a di Palermo  (Co.R.I.).}\\
\footnotesize{Dipartimento di Metodi e Modelli Matematici}\vspace*{-1,5mm}\\
\footnotesize{Viale delle Scienze Ed. 8, I-90128 Palermo (Italy)}\vspace*{-1,5mm}\\
\footnotesize{gfalcone@unipa.it}}

\vskip1cm\begin{abstract} \noindent This paper is the extended abstract of a talk given at the conference \em Combinatorics 2008\em , Costermano (Italy) 22-28 June 2008. We try to embed a $t$-design $\cal D$ in a finite commutative group in such a way that the sum
of the $k$ points of a block is zero. We can compute the number of blocks of the \em boolean \em $2$-design having all the non zero vectors of ${\Bbb Z}_2^n$ as the set of points and the $k$-subsets of elements the sum of which is zero as blocks.

\end{abstract}

\Large\vskip1cm\noindent
{\underline{\bf 1. Preliminaries.}}

\large\bigskip\noindent
We recall that a $t-(v,k,r_t)$ design ${\cal D}=({\cal P},{\cal B})$ is a finite set $\cal P$ with $|{\cal P}|=v$, the elements of which are called \em points\em , together with a family ${\cal B}$ of subsets of $\cal P$, called \em blocks\em , such that any block contains exactly $k$ points and $t$ points are contained in exactly $r_t$ blocks. For any $s< t$, any $t-(v,k,r_t)$ design ${\cal D}$ is an $s-(v,k,r_s)$ design with
$$r_s=r_t\;{{(v-s)(v-s-1)\cdots (v-t+1)}\over{(k-s)(k-s-1)\cdots (k-t+1)}}.$$
The number $r_1$ of blocks containing any given point is consequently a constant. It is usual to denote $r_1$ simply by $r$ and $r_2$ by $\lambda$.
If we denote by $b=|{\cal B}|$ the number of blocks, a necessary condition for the existence of a $t-(v,k,r_t)$ design is $vr=bk$.

\smallskip\noindent
Labelling the points of $\cal P$ with $P_1,\cdots ,P_v$ and the blocks of $\cal B$ with ${\mathfrak{b}}_1,\cdots ,{\mathfrak{b}}_b$, the \em incidence matrix \em $A=(a_{ij})$ of $\cal D$ is defined putting $$a_{ij}=\left\{\begin{array}{cc}
                                              1 & \mbox{ if } P_i\in {\mathfrak{b}}_j \\
                                              0 & \mbox{ if } P_i\notin {\mathfrak{b}}_j
                                            \end{array}.
\right.$$
It holds $A^{\sf T}A=(r-\lambda)I+\lambda J$, where $J$ is the $v\times v$ matrix the entries of which are all equal $1$. It follows that $|A^{\sf T}A|=rk(r-\lambda)^{(v-1)}$. A $2$-design where $v=b$, or equivalently $r=k$, is called \em symmetric \em and its incidence matrix is such that
$AA^{\sf T}=A^{\sf T}A$.

\smallskip\noindent
The \em complementary design \em of ${\cal D}=({\cal P},{\cal B})$ is the $t-(v,v-k,\tilde{r_t})$ design $\tilde{{\cal D}}=({\cal P},\tilde{{\cal B}})$ where $\tilde{{\cal B}}$ is the set of $(v-k)$-tuples of $\cal P$ which are the complement ${\cal P}\setminus {\mathfrak{b}}$ of a block ${\mathfrak{b}}\in\cal B$. Consequently, we have $\tilde{r_t}=(b-2r+\lambda){{(k-2)\cdots (k-t+1)}\over{(v-2)\cdots (v-t+1)}}$.

\smallskip\noindent
The \em supplementary design \em of ${\cal D}=({\cal P},{\cal B})$ is the $t-(v,k,{v-t\choose k-t} - r_t)$ design $\breve{{\cal D}}=({\cal P},\breve{{\cal B}})$ where $\breve{{\cal B}}$ is the set of unordered $k$-tuples of distinct points of $\cal P$ which are not blocks of $\cal B$.

\smallskip\noindent
{\bf 1. Remark:} The $t-(v,k,r_t)$ design ${{\cal D}}$ is the
complementary design of $\tilde{{\cal D}}$ as well as the
supplementary design of $\breve{{\cal D}}$. The complementary
design of the supplementary design of ${{\cal D}}$ is a
$t-(v,v-k,\breve{r}_t)$ design which is equal to the supplementary
design of the complementary design of ${{\cal D}}$.\hfill$\Box$

\smallskip\noindent
The \em derived design of ${\cal D}=({\cal P},{\cal B})$ at the point $P$ \em is the design ${\sf Der}_P{\cal D}=({\cal P}\setminus P,{\sf Der}_P{{\cal B}})$ where ${\sf Der}_P{{\cal B}}=\{{\mathfrak b}\setminus P:\;P\in{\mathfrak b}\in{\cal B}\}$.

\smallskip\noindent
Lastly, we recall that a \em Steiner $k$-tuple system \em is a $t-(v,k,r_t)$ design with $t=k-1$ and $r_t=1$. Among the Steiner quadruple systems, i. e. $3-(v,4,1)$ designs, we find the \em boolean quadruple system of order $2^n$\em , which is defined, for $n\geq 3$, as the $3-(2^n,4,1)$ design obtained putting ${\cal B}$ to be the set of all quadruples of distinct vectors of ${\cal P}={\Bbb Z}_2^n$ the sum of which is zero.

\Large\vskip1cm\noindent {\underline{\bf 2. Embedding in a
group.}}

\large\bigskip\noindent
We start by asking a natural question: what designs are subsets of a finite commutative group, so that the sum of the elements in a block is a constant? If this constant is zero, there is a unique way to define such a group.
Let $\mathfrak{G}$ be the free commutative group generated by the $v$ points of $\cal P$ and let $\mathfrak{R}$ be the subgroup of $\mathfrak{G}$ generated by the $b$ elements of the form
$$\sum_{X\in {\mathfrak{b}}_j}X \quad (j=1,\cdots ,b),$$
where ${\mathfrak{b}}_j$ is a block of $\cal B$. The subgroup
$\mathfrak{R}$ is clearly generated by the rows of the incidence
matrix of $\cal D$. Finally, define the group
${\mathfrak{G}}_{\cal D}=\mathfrak{G}/\mathfrak{R}$ and consider
the map ${\cal P}\longrightarrow {\mathfrak{G}}_{\cal D}$,
$X\mapsto x=X+\mathfrak{R}$.

\bigskip\noindent{\textbf{2. Example:}}
Let $\cal D$ be the $2-(9,3,1)$ Steiner triple system. After reducing the rows of the incidence matrix by elementary integer linear combinations
$${\small{A=\left(
    \begin{array}{ccccccccc}
      1 & 1 & 1 & 0 & 0 & 0 & 0 & 0 & 0 \\
      0 & 0 & 0 & 1 & 1 & 1 & 0 & 0 & 0 \\
      0 & 0 & 0 & 0 & 0 & 0 & 1 & 1 & 1 \\
      1 & 0 & 0 & 1 & 0 & 0 & 1 & 0 & 0 \\
      0 & 1 & 0 & 0 & 1 & 0 & 0 & 1 & 0 \\
      0 & 0 & 1 & 0 & 0 & 1 & 0 & 0 & 1 \\
      1 & 0 & 0 & 0 & 0 & 1 & 0 & 1 & 0 \\
      0 & 1 & 0 & 1 & 0 & 0 & 0 & 0 & 1 \\
      0 & 0 & 1 & 0 & 1 & 0 & 1 & 0 & 0 \\
      1 & 0 & 0 & 0 & 1 & 0 & 0 & 0 & 1 \\
      0 & 1 & 0 & 0 & 0 & 1 & 1 & 0 & 0 \\
      0 & 0 & 1 & 1 & 0 & 0 & 0 & 1 & 0 \\
    \end{array}
  \right)\rightsquigarrow\left(
                           \begin{array}{ccccccccc}
                             1 & 0 & 0 & 0 & 0 & 1 & 0 & 1 & 0 \\
                             0 & 1 & 0 & 0 & 0 & -2 & 0 & -1 & -1 \\
                             0 & 0 & 1 & 0 & 0 & 1 & 0 & 0 & -2 \\
                             0 & 0 & 0 & 1 & 0 & -1 & 0 & -2 & -1 \\
                             0 & 0 & 0 & 0 & 1 & 2 & 0 & 2 & 1 \\
                             0 & 0 & 0 & 0 & 0 & 3 & 0 & 3 & 0 \\
                             0 & 0 & 0 & 0 & 0 & 0 & 1 & 1 & -2 \\
                             0 & 0 & 0 & 0 & 0 & 0 & 0 & 3 & 0 \\
                             0 & 0 & 0 & 0 & 0 & 0 & 0 & 0 & 3 \\
                             0 & 0 & 0 & 0 & 0 & 0 & 0 & 0 & 0 \\
                             0 & 0 & 0 & 0 & 0 & 0 & 0 & 0 & 0 \\
                             0 & 0 & 0 & 0 & 0 & 0 & 0 & 0 & 0 \\
                           \end{array}
                         \right),}}
$$
we find that ${\mathfrak{G}}_{\cal D}=({\Bbb Z}_3)^3$, and the
(images of the) points of $\cal D$ are:
$$\begin{array}{ccc}
                 P_1=(0,2,2) & P_2=(1,1,2) & P_3=(2,0,2) \\
                 P_4=(1,2,1) & P_5=(2,1,1) & P_6=(0,0,1) \\
                 P_7=(2,2,0) & P_8=(0,1,0) & P_9=(1,0,0).
               \end{array}
$$
Thus ${\cal D}=\{(x_1,x_2,x_3)\in ({\Bbb Z}_3)^3:\;x_1+x_2+x_3=1\}$. This gives evidence of the fact that $\cal D$ is isomorphic to the
point-line design of the affine plane of order $3$. Now we consider the $2-(9,6,5)$ design of pairs of parallel lines in the affine plane of order $3$. Here, after reducing the rows of the incidence matrix
$${\small{A=\left(
    \begin{array}{ccccccccc}
      1 & 1 & 1 & 1 & 1 & 1 & 0 & 0 & 0 \\
      1 & 1 & 1 & 0 & 0 & 0 & 1 & 1 & 1 \\
      1 & 1 & 0 & 1 & 1 & 0 & 1 & 1 & 0 \\
      1 & 0 & 1 & 1 & 0 & 1 & 1 & 0 & 1 \\
      1 & 1 & 0 & 0 & 1 & 1 & 1 & 0 & 1 \\
      1 & 0 & 1 & 1 & 1 & 0 & 0 & 1 & 1 \\
      1 & 1 & 0 & 1 & 0 & 1 & 0 & 1 & 1 \\
      1 & 0 & 1 & 0 & 1 & 1 & 1 & 1 & 0 \\
      0 & 1 & 1 & 0 & 1 & 1 & 0 & 1 & 1 \\
      0 & 0 & 0 & 1 & 1 & 1 & 1 & 1 & 1 \\
      0 & 1 & 1 & 1 & 1 & 0 & 1 & 0 & 1 \\
      0 & 1 & 1 & 1 & 0 & 1 & 1 & 1 & 0 \\
    \end{array}
  \right)\rightsquigarrow\left(
                           \begin{array}{ccccccccc}
                             1 & 0 & 0 & 0 & 0 & 1 & 0 & -2 & 0 \\
                             0 & 1 & 0 & 0 & 0 & -2 & 0 & 2 & -1 \\
                             0 & 0 & 1 & 0 & 0 & 1 & 0 & 0 & -2 \\
                             0 & 0 & 0 & 1 & 0 & -1 & 0 & 1 & -1 \\
                             0 & 0 & 0 & 0 & 1 & 2 & 0 & -1 & -2 \\
                             0 & 0 & 0 & 0 & 0 & 3 & 0 & -3 & 0 \\
                             0 & 0 & 0 & 0 & 0 & 0 & 1 & 1 & -2 \\
                             0 & 0 & 0 & 0 & 0 & 0 & 0 & 3 & -3 \\
                             0 & 0 & 0 & 0 & 0 & 0 & 0 & 0 & 6 \\
                             0 & 0 & 0 & 0 & 0 & 0 & 0 & 0 & 0 \\
                             0 & 0 & 0 & 0 & 0 & 0 & 0 & 0 & 0 \\
                             0 & 0 & 0 & 0 & 0 & 0 & 0 & 0 & 0 \\
                           \end{array}
                         \right)
}}$$ we find that ${\mathfrak{G}}_{\cal D}={\Bbb Z}_2\oplus({\Bbb
Z}_3)^3$, and the points of $\cal D$ are:
$$\begin{array}{ccc}
                 P_1=(1;0,2,2) & P_2=(1;1,1,2) & P_3=(1;2,0,2) \\
                 P_4=(1;1,2,1) & P_5=(1;2,1,1) & P_6=(1;0,0,1) \\
                 P_7=(1;2,2,0) & P_8=(1;0,1,0) & P_9=(1;1,0,0)
               \end{array}
$$
We remark that it is unexpected that the coordinates of the points
are precisely the same as before, apart from the first one.
Moreover, even if the first coordinate is constant, it seems to
plays an important r\^{o}le, since it distinguish blocks from the
first case: note that $(P_1+P_2+P_3)$ and $(P_4+P_5+P_6)$ here are
not zero, whereas
$(P_1+P_2+P_3)+(P_4+P_5+P_6)=(0;0,0,0)$.\hfill$\Box$

\bigskip\noindent
In the following proposition, which has clearly a connection with
the computation of the $p$-rank of the incidence matrix in
\cite{Hamada} and \cite{DHV}, shows that the exponent of
${\mathfrak{G}}_{\cal D}$ divides $k(r-\lambda)$.

\bigskip\noindent{\textbf{3. Proposition:}}
\em For any $x=X+\mathfrak{R}\in {\mathfrak{G}}_{\cal D}$ we have
$k(r-\lambda)x=0$. If ${\cal D}$ has a partition in blocks, then
$(r-\lambda)x=0$.\em

\smallskip\noindent
{\em Proof.} Summing the points of the $r$ blocks through any
given point $X$ we find
$$(r-\lambda)x+\lambda\cdot\sum_{y=Y+{\mathfrak R}}y=0,$$
hence $(r-\lambda)x=0$ , if ${\cal D}$ has a partition in blocks.
Otherwise, let $X_1,\cdots ,X_k$ be the points of a given block,
then
$$0=(r-\lambda)\sum_{x_i=X_i+{\mathfrak R}} x_i=-\, k\cdot \lambda\cdot\sum_{y=Y+{\mathfrak
R}}y.$$ As $k\cdot \lambda\cdot\sum_{y=Y+{\mathfrak R}}y=0$, the
assertion follows.\hfill$\Box$

\bigskip\noindent
We remark that the injectivity of the map ${\cal P}\longrightarrow
{\mathfrak{G}}_{\cal D}$, $X\mapsto x=X+\mathfrak{R}$, is not
always guaranteed. In particular we find that, if $v\equiv
1\quad(12)$, any Steiner triple system $\cal D$ of cardinality $v$
is not embeddable in ${\mathfrak{G}}_{\cal D}$. In fact, as a
consequence of \cite{DHV}, we have that the $b\times v$ incidence
matrix of $\cal D$ has rank $v$ over any field of characteristic
$p\not=3$ and it has rank $v-1$ over a field of characteristic
$p=3$. This forces ${\mathfrak{G}}_{\cal D}$ to be a cyclic group
of a $3$-power order, but by the above Proposition the exponent of
${\mathfrak{G}}_{\cal D}$ must be a divisor of $3$. Since $\cal D$
has more than three points, we can see that $\cal P$ is not
embeddable in ${\mathfrak{G}}_{\cal D}$.

\bigskip\noindent
The following proposition, though trivial, has a fine corollary.

\bigskip\noindent{\textbf{4. Proposition: }}
\em The map is injective if and only if for any
${\mathbf{v}}\in{\Bbb Z}^b$ and for any permutation matrix $H$, we
have\em  $${\mathbf{v}}A\not=(1,-1,0,\cdots ,0)H.$$
 \hfill$\Box$

\bigskip\noindent{\textbf{5. Corollary: }}(Irrationality condition for the injectivity):
\em If the map is not injective, then there exists a vector
${\mathbf{v}}\in{\Bbb Z}^b$
such that\em
$$\langle{\mathbf{v}},{\mathbf{v}}\rangle={\mathbf{v}}AA^{\sf
T}{\mathbf{v}}^{\sf T}=2.$$

\smallskip\noindent
{\em Proof.} The assertion follows from the fact that, for any
permutation matrix $H$, we have $HH^{\sf T}=I$. \hfill$\Box$

\bigskip\noindent{\textbf{6. Corollary: }}
\em A symmetric $2-(v,k,\lambda)$ design ${\cal D}$ is embeddable
in ${\mathfrak{G}}_{\cal D}$ (unless $\cal D$ is a triangle).\em

\smallskip\noindent
{\em Proof.} It is known that a $2$-design is symmetric if and
only if $AA^{\sf T}=A^{\sf T}A=(k-\lambda)I+\lambda J$. Hence, for
any non-zero integer $b$-tuple ${\mathbf{v}}=(x_1,\cdots ,x_b)$,
we have
$$\langle{\mathbf{v}},{\mathbf{v}}\rangle={\mathbf{v}}AA^{\sf T}{\mathbf{v}}^{\sf T}=(k-\lambda)\sum_{i=1}^bx_i^2+\lambda\left(\sum_{i=1}^bx_i\right)^2>2.$$
 \hfill$\Box$

\bigskip\noindent
As pointed out in Proposition 3, it is necessary that the exponent
of ${\mathfrak{G}}_{\cal D}$ is a divisor of $k(r-\lambda)$. Now
we show that for the existence of a $2-(v,k,\lambda)$ design
having a $p$-group as the set of points it is sufficient that $kx=0$.

\smallskip\noindent
{\bf 7. Proposition:} \em Let $\cal P$ be the Galois field with
$q=p^t$ elements, $p\geq 2$ a prime. For any $k=mp$, with $2<k<q$,
let ${\cal B}$ be the family of unordered $k$-tuples of distinct
elements of $\cal P$ the sum of which is zero. Then ${\cal D}=({
\cal{P}},{ \cal{B}})$ is a $2-(q,k,\lambda)$ design. \em

\smallskip\noindent
{\em Proof.} As the sum of the elements in the ground field is
zero, we have that ${\cal B}\not=\emptyset$. Let now
$P_1,P_2\in{\mathfrak b}\in {\cal B}$, let $Q_1,Q_2\in {\cal P}$
and let $\rho(X)={A}X+T$, with $A, T\in {\cal P}$ and $A\not= 0$,
be the affinity of the affine line defined on ${\cal P}$, such
that $\rho(P_i)=Q_i$. Then
$$\sum_{{{Q}}\in\rho({\mathfrak{b}})}Q=\sum_{P\in {\mathfrak{b}}}\rho(P)=\sum_{P\in {\mathfrak{b}}}({A}P+T)=k\, T.$$
As $k\equiv 0\;(p)$ we find that $\rho({\mathfrak{b}})$ is in
$\cal B$, hence the number of unordered $k$-tuples of $\cal B$
containing $Q_1,Q_2$ is equal to the number of unordered
$k$-tuples of $\cal B$ containing $P_1,P_2$.\hfill$\Box$

\smallskip\noindent
{\bf 8. Remark:} A remarkable case is when $k=p$ and ${\cal P}$ is
an $n$--dimensional affine space over a Galois field with $p^m$
elements, $m>1$.\hfill$\Box$

\Large\vskip1cm\noindent
{\underline{\bf 3. Boolean designs.}}

\large\bigskip\noindent In this section we consider the case
$p=2$.

\bigskip\noindent
{\bf 9. Proposition:} \em Let $\cal P$ be a $n$--dimensional
affine space over the Galois field with $2$ elements. For any
$k=2m$, with $2<k<2^n$, let ${\cal B}$ be the family of unordered
$k$-tuples of distinct elements of $\cal P$ the sum of which is
zero. Then ${\cal D}=({ \cal{P}},{ \cal{B}})$ is a $3-(2^n,k,r_3)$
design. \em

\smallskip\noindent
{\em Proof.} Since a line in $\cal P$ cannot contain three points
and the group of affinity is transitive on the triangles, we can
move any three distinct points of $\cal P$ onto any three distinct
points of $\cal P$ with an affinity. The proof follows as in
Proposition 7.\hfill$\Box$

\smallskip\noindent
In addition to the above one, in the following proposition we
apply the fact that $GL_n(2)$ is transitive on the pairs of
non-zero distinct vectors of ${\Bbb Z}_2^n$. For $k=3$ we find
precisely the classical point-line design of a projective space on
${\Bbb Z}_2^n$.

\smallskip\noindent
{\bf 10. Proposition:} \em Let ${\cal P}$ be the set of non--zero
vectors of ${\Bbb Z}_2^n$ and, for any $k=2,3,\cdots,2^n-2 $, let
${\cal B}_k$ be the family of unordered $k$-tuples of distinct
elements of $\cal P$ the sum of which is zero. Then ${\cal
D}_k=({\cal P},{\cal B}_k)$ is a $2-(v=2^n-1,k,\lambda_k)$ design
such that $(k+1)b_{k+1}+b_k+(v-k+1)b_{k-1}={{v}\choose k}$.
Consequently we have $b_k={{v}\choose k}\alpha_{\lfloor{{k-1}\over
2}\rfloor}$ where
$$\begin{array}{ll}
             \alpha_h & ={1\over{v-2h}}\left(1-\sum_{i=0}^{h-2}(-1)^i\prod_{j=0}^i{{1+2(h-j)}\over{v-2(h-j-1)}}\right) \\
             \\
              & ={1\over{v-2h}}\sum_{i=0}^{h-1}(-1)^i{{(2h+1)!!}\over{\big(2(h-i)+1\big)!!}}\cdot {{(v-2h)!!}\over{\big(v-2(h-i-1)\big)!!}}.
           \end{array}
$$
\em \hfill$\Box$

\smallskip\noindent
{\bf 11. Remark:}
 According to \cite{Schatz} (see also sequences A010085-89 in \cite{OEIS}), the numbers $b_k$ of blocks of ${\cal D}_k$ are equal to the weights of the \em huge \em $(2^n-1, 2^n-n-1,3)$--Hamming code $C$, whereas the numbers $\breve{b}_k$ of blocks of $\breve{{\cal D}}_k$ are the numbers of weights $k$ vectors of a $C$ which belong to weight $1$ cosets of $C$.
 \hfill$\Box$

\smallskip\noindent
{\bf 12. Remark:}
 The supplementary design $\breve{{\cal D}}_k$ of the above one is the $2-(2^n-1,k,\breve{r}_k)$ design defined, for any $1<k<2^n-1$, on the set ${\cal P}$ of non--zero vectors of ${\Bbb Z}_2^n$ by the family ${\cal B}$ of unordered $k$-tuples of distinct elements of $\cal P$ the sum of which is different from zero.
 \hfill$\Box$

\smallskip\noindent
{\bf 13. Remark:} For $k=4$ in Proposition 9 we get the boolean
system of order $2^n$. We remark that in the affine space ${\Bbb
Z}_2 ^n$ a necessary and sufficient condition for four distinct
points to lie in an affine plane is that their sum is zero, that
is the boolean quadruple system of order $2^n$ is the \em
classical \em design of two-dimensional subspaces of an affine
space over ${\Bbb Z}_2$. Consider now the design where $k=8$. We
remark that, in the affine spaces ${\Bbb Z}_2^4$ and ${\Bbb
Z}_2^5$, a necessary and sufficient condition for eight distinct
points $P_1,\cdots ,P_8$ to lye in a two disjoint planes is that
their sum is zero. To see this, we can assume that any four of
${\cal I}=\{P_1,\cdots P_8\}$ do not lie in a plane. Taking the
sum of any three points of $\cal I$, we get ${2^n\choose 3}$
further points, not in $\cal I$, which are mutually distinct,
otherwise two of the eight points of $\cal I$ are equal. For
$n=4,5$ we get then a contradiction. Hence the design we get for
$n=4,5$ and $k=8$ is the \em classical \em design of disjoint
pairs of two-dimensional subspaces of an affine space over ${\Bbb
Z}_2$. Things change for ${\Bbb Z}_2^6$, because the following
$8$-ple, the sum of which is zero,

$(0,0,0,0,0,0)$\quad $(1,0,0,0,0,0)$\quad $(0,1,0,0,0,0)$\quad $(0,0,1,0,0,0)$

$(0,0,0,1,0,0)$\quad $(0,0,0,0,1,0)$\quad $(0,0,0,0,0,1)$\quad $(1,1,1,1,1,1)$,

\smallskip\noindent
is not the disjoint union of two affine subplanes.\hfill$\Box$

\bigskip\noindent {\bf 14. Definition:} \em Let ${\cal D}_k=({\cal
P},{\cal B}_k)$ be the $2-(2^n-1,k,\lambda_k)$ design where ${\cal
P}$ is the set of non--zero vectors of ${\Bbb Z}_2^n$ and ${\cal
B}_k$ is the family of unordered $k$-tuples of distinct elements
of $\cal P$ the sum of which is zero. We say that the block
${\mathfrak{b}}\in {\cal B}_k$ is \em reducible \em if it is the
union of two disjoint blocks ${\mathfrak{b}}_1\in{\cal B}_{k_1}$,
${\mathfrak{b}}_2\in{\cal B}_{k_2}$, where $k_1+k_2=k$. \em

\bigskip\noindent{\textbf{15. Proposition:}}
\em Let $\{{\mathbf{e}}_i:\, i=1,2,\cdots , n\}$ be the canonical basis of ${\Bbb Z}_2^n$ and let ${\mathfrak{c}}_k=\{{\mathbf{e}}_1,\cdots ,{\mathbf{e}}_{k-1},{\mathbf{e}}_1+{\mathbf{e}}_2+\cdots +{\mathbf{e}}_{k-1}\}$.
Then any block ${\mathfrak{b}}\in{\cal B}_k$ is irreducible if and only if ${\mathfrak{b}}$ is contained in the orbit of ${\mathfrak{c}}_k$ under $GL_n(2)$.
\em

\smallskip\noindent
{\em Proof.} Let ${\mathfrak{b}}=\{P_1,\cdots ,P_{k-1},P_k=\sum_{j=1}^{k-1}P_j\}$ be an irreducible block. The claim follows if we show that
$\{P_1,\cdots ,P_{k-1}\}$ are linearly independent. By contradiction, we may assume without loss of generality that $P_{k-1}=\sum_{j=1}^{k-2}\alpha_jP_j$, with $\alpha_j=0,1$, but not all zero. If any $\alpha_j=1$, then $P_k={\mathbf{0}}$, a contradiction. Thus some $\alpha_j=0$. But this forces ${\mathfrak{b}}$ to be reducible, a contradiction.\hfill$\Box$

\smallskip\noindent
{\bf 16. Corollary:} \em  The number of irreducible blocks
${\mathfrak{b}}\in{\cal B}_k$ is $\prod_{i=1}^{k}
(2^n-2^{i-1})$.\em \hfill$\Box$

\smallskip\noindent
{\bf 17. Remark:}
 The family ${\cal B}$ of unordered $k$-tuples of linearly dependent vectors of $\cal P$ is such that ${\mathfrak{b}}\in{\mathcal{B}}$ if and only
 if ${\mathfrak{b}}$ contains a $h$-tuple of non-zero vectors the sum of which is zero, for some $h\leq k$.
 This shows that ${\cal D}=({\cal P},{\cal B})$ is a $2$-design,
 the supplementary of which is the $2$-design defined by the family $\breve{{\cal B}}$ of unordered $k$-tuples of
 linearly independent vectors of $\cal P$.
\hfill$\Box$

\end{document}